# Computing Roots of Directed Graphs is Graph Isomorphism Hard

Martin Kutz[*]
Freie Universität Berlin
kutz@math.fu-berlin.de

October 29, 2018

**Abstract**

The $k$-th power $D^k$ of a directed graph $D$ is defined to be the directed graph on the vertices of $D$ with an arc from $a$ to $b$ in $D^k$ iff one can get from $a$ to $b$ in $D$ with exactly $k$ steps. This notion is equivalent to the $k$-fold composition of binary relations or $k$-th powers of Boolean matrices.

A $k$-th root of a directed graph $D$ is another directed graph $R$ with $R^k = D$. We show that for each $k \geq 2$, computing a $k$-th root of a directed graph is at least as hard as the graph isomorphism problem.

**Keywords.** directed graph, graph power, root, binary relation, Boolean matrix, graph isomorphism, computational complexity

**AMS classification.** 05C12, 05C20, 05C60, 68Q17, 05C50, 15A23, 06E99

## 1 Introduction

Let $D$ be a directed graph without multiple arcs, it may have loops. We define the $k$-th power $D^k$ of $D$, $k \in \mathbb{N}$, to be the directed graph on the same vertex set and with an arc from $a$ to $b$ in $D^k$ if and only if one can get from $a$ to $b$ in $D$ with *exactly* $k$-steps, possibly visiting some vertices several times. Compare Figure 1.

This definition arises naturally from the connection between directed graphs, binary relations, and Boolean matrices. Interpreted as a binary relation, the $k$-th power $D^k$ is simply the $k$-fold composition of the relation $D$.

---

[*]member of the European graduate school "Combinatorics, Geometry, and Computation" supported by the Deutsche Forschungsgemeinschaft, grant GRK 588/1

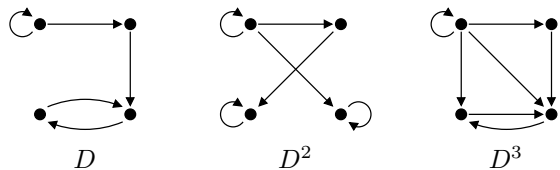

Figure 1: Powers of a digraph.

And if we compute Boolean powers of the zero-one adjacency matrix $A$ of a digraph $D$, that is, we perform ordinary matrix multiplication with $+$ and $\cdot$ replaced by the Boolean operations $\vee$ and $\wedge$, we see that the adjacency matrix of $D^k$ is just $A^k$. So powers of digraphs, binary relations, and Boolean matrices are equivalent concepts [RW88].

A $k$-th root of a digraph $D$ is a digraph $R$ satisfying $R^k = D$. We investigate the computational complexity of computing roots of directed graphs.

In the open problems section of his book [Kim82], Kim asks whether for the special case $k = 2$, such roots can be computed in polynomial time or whether this problem it is perhaps $\mathcal{NP}$-complete. (Actually, he poses this question in terms of square Boolean matrices, not digraphs.) We give a partial answer to this question by relating digraph powers to graph isomorphism. A computational problem that is at least as hard as finding isomorphisms between graphs is called *isomorphism hard*. (We will cover the precise definition of these complexity notions later.) Here is our main result.

**Theorem 1.** *For each single parameter $k \geq 2$, the problem of deciding whether a given directed graph has a $k$-th root is isomorphism hard.*

Graph isomorphism is a famous candidate for a problem strictly between $\mathcal{P}$ and $\mathcal{NP}$-completeness, provided that actually $\mathcal{P} \neq \mathcal{NP}$. Problems of the same complexity as graph isomorphism are called *isomorphism complete*. So we are left with the natural question what the precise complexity of computing digraph roots might be. Whether it is isomorphism complete, $\mathcal{NP}$-complete, or maybe of intermediate complexity.

As a first small step towards a possible isomorphism completeness proof of the $k$-th root problem, we show that for a certain class of digraphs that play an important role in the proof of Theorem 1, root finding is actually isomorphism complete. Maybe this result can be extended to larger classes of digraphs that maintain this property.

Our notion of digraph powers must not be confused with a possibly more common one considered in, for example, [MS94] and [Fle74]. In [MS94], Motwani and Sudan prove the problem of deciding whether an *undirected* graph has a square root $\mathcal{NP}$-complete. But they define the square of a graph through paths of length *at most* two. In that setting, edges cannot vanish through squaring, as happened in Figure 1. In other words, they consider square roots in the class of symmetric and reflexive relations. It seems that this yields an essentially different problem.

## 2  Basic Concepts

Let us begin with a precise definition of our notion of directed graphs. A *directed graph* $D$ (*digraph*, for short) is a finite set $V(D)$ of *vertices* together with a set $A(D) \subseteq V(D) \times V(D)$, the *arcs* of $D$. Note that this definition excludes multiple arcs but allows loops, as intended. For brevity, we simply write $ab$ for pairs $(a, b) \in V(D) \times V(D)$. The expression $a \to b$ shall indicate that $ab \in A(D)$ wherever the reference to some digraph $D$ is clear.

A *walk* of length $k$ in a digraph $D$ is a sequence $(a_0, a_1, \ldots, a_k)$ of vertices with $a_i a_{i+1} \in A(D)$ for $0 \leq i < k$. Note that a vertex may appear several times in a walk. A walk is a *path* if all its vertices are different. A walk is a *cycle* if $a_0 = a_k$. It may be instructive to restate the definition of digraph powers in terms of walks. For each $k \in \mathbb{N}$, the $k$-th power $D^k$ of $D$ is the digraph on the same vertices as $D$ and with an arc from $a$ to $b$ in $D^k$ if and only if there is a walk of length *exactly* $k$ from $a$ to $b$ in $D$.

We fix some further notation for digraphs. Let $D$ be some digraph. We let $O_D(a) := \{x \in X \mid ax \in A(D)\}$ denote the set of *out-neighbors* and $I_D(a) := \{x \in X \mid xa \in A(D)\}$ the set of *in-neighbors* of $a \in V(D)$. Their cardinalities are denoted by $\delta_D^+(a) := |O_D(a)|$ and $\delta_D^-(a) := |I_D(a)|$. If the reference to the digraph in question is clear, we may also omit the subscript. For subsets $U \subseteq V(D)$, we let $O_D(U) := \bigcup_{u \in U} O_D(u)$ and similarly for $I_D$. Note that our definition of $I$ and $O$ differs from the usual convention in that our $U$-neighborhoods may contain elements from $U$. We shall often need iterated neighbourhoods and hence define recursively $O^k(U) := O(O^{k-1}(U))$ and $I^k(U) := I(I^{k-1}(U))$ with $O^1(U) = O(U)$ and $I^1(U) = I(U)$. For convenience, we also let $O^0(U) = I^0(U) = U$.

If again $U$ is a subset of $V(D)$, the expression $D[U]$ denotes the *induced* subgraph on $U$. That is, the digraph with vertex set $U$ and with $ab \in A(D[U])$ if and only if $ab \in A(D)$.

A digraph $D$ is *weakly connected* if for any two vertices $a, b \in V(D)$ there exists a sequence $a = x_0, x_1, \ldots, x_m = b$ of vertices so that $x_i x_{i+1} \in A(D)$ or $x_{i+1} x_i \in A(D)$ for $0 \leq i < m$. In other words, we can get from $a$ to $b$ using arcs in any direction. A *weakly connected component* of a digraph $D$ is a maximal weakly connected induced subgraph of $D$.

## Graph Isomorphism

An isomorphism between two digraphs $D_1$ and $D_2$ is a bijection $\varphi \colon D_1 \to D_2$ so that $ab \in A(D_1)$ iff $\varphi(a)\varphi(b) \in A(D_2)$. Two digraphs are isomorphic if there exists an isomorphism between them.

Given two directed graphs, the graph isomorphism problem asks whether these graphs are isomorphic or not. We remark that usually, undirected graphs are considered, but a few simple transformations show that with respect to their computational complexity, the two problems are equivalent [KST93]. So we use directed graphs because they are better suited for our needs.

Since there exist (slightly) different concepts for comparing the complexity of computational problems, we briefly recall the one we shall use. A *Karp reduction* (or *many-one reduction*) from a problem $X$ to another problem $Y$ is a polynomial-time computable mapping from the instances of $X$ to those of $Y$ so that positive instances map to positive instances and negative to negative ones. If such a mapping exists, we say that problem $X$ *reduces* to problem $Y$. It means that solving $X$ is at most as hard as solving $Y$. A problem that graph isomorphism reduces to is called *isomorphism hard*. If such a problem also reduces to graph isomorphism, we call it *isomorphism complete*.

Several computational problems have been shown isomorphism complete, most of these are isomorphism problems for other algebraic or combinatorial structures. For example, isomorphism of semigroups and finite automata [Boo78], convex polytopes [KS01], or finitely represented algebras. Other problems ask for properties of the automorphism group of a graph, for example, computing the size of the automorphism group or its orbits [Mat79]. (We note that the latter two problems are known to be graph isomorphism complete only in the weaker sense of Turing reductions, which are defined through oracle machines.) Finally, several restrictions of the graph isomorphism problem are known to remain isomorphism complete. Among these are regular graphs [Mat79], line graphs, and bipartite graphs. For deeper informations about the complexity of graph isomorphism we refer to the book [KST93].

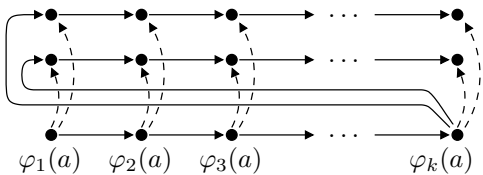

Figure 2: Constructing a $k$-th root (continuous lines) for a disjoint union of $k$ isomorphic digraphs (dashed lines).

## 3 The Reduction

For the proof of Theorem 1 we need a Karp reduction from graph isomorphism to $k$-th roots of digraphs. We have to transform a pair $D_1, D_2$ of digraphs into a single digraph $E$ so that $D_1$ and $D_2$ are isomorphic if and only if $E$ has a $k$-th root. And we need this mapping for each parameter $k \geq 2$. Our search for such a transformation will be guided by the following observation about disjoint unions of isomorphic digraphs.

**Proposition 1.** *Let $D = D_1 \dot\cup D_2 \dot\cup \cdots \dot\cup D_k$ be the disjoint union of $k$ isomorphic digraphs $D_1, \ldots, D_k$. Then $D$ has a $k$-th root.*

*Proof.* We construct a digraph $R$ on the vertices of $D$ with $R^k = D$. Let $D_0$ be a further isomorphic copy of some of the $D_i$. There exist digraph isomorphisms
$$\varphi_i \colon D_0 \to D_i, \quad 1 \leq i \leq k.$$
For each vertex $a$ of $D_0$, we let $R$ contain the path
$$\varphi_1(a) \to \varphi_2(a) \to \cdots \to \varphi_k(a)$$
and additionally the arcs
$$\varphi_k(a) \to \varphi_1(b) \quad \text{for all } b \in O_{D_0}(a). \tag{1}$$
Figure 2 shows a local picture of this construction.

We claim that $R^k = D$. To see this, pick any $v \in D_i$, $1 \leq i \leq k$, and let $a := \varphi_i^{-1}(v)$. We have
$$\begin{aligned}
O_R^k(v) &= O_R^{i-1}\big(O_R(O_R^{k-i}(v))\big) \\
&= O_R^{i-1}\big(O_R(\varphi_k(a))\big) \\
&= O_R^{i-1}\big(\varphi_1(O_{D_0}(a))\big) \quad \text{by (1)} \\
&= \varphi_i(O_{D_0}(a)) \\
&= O_{D_i}(v) = O_D(v).
\end{aligned}$$
□

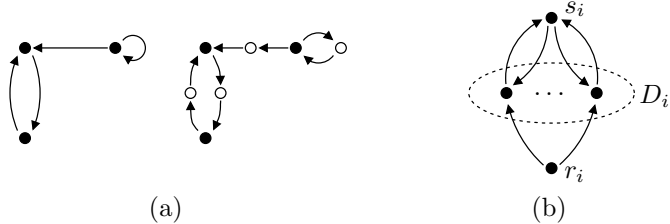

(a)                (b)

Figure 3: (a) A digraph and its complete subdivision. The black vertices of the subdivision digraph form a core. (b) The suspension of a digraph $D_i$.

Proposition 1 indicates how to reduce graph isomorphism to digraph roots. Given two digraphs $D_1, D_2$, we form the disjoint union $D$ of $D_1, D_2$, and $k-2$ copies of $D_2$. If $D_1$ and $D_2$ are isomorphic then $D$ has a $k$-th root. What if they are not. Can $D$ still have a $k$-th root then? Unfortunately, the answer is "yes". One easily finds such examples when one of $D_1$ and $D_2$ contains isolated vertices. So the converse direction of our construction does not work yet. But we shall be able to fix this. All we need is some additional structure in the digraph $D$. In order to describe the reduction in full detail, we must first introduce a few concepts.

### Subdivisions, Cores, and Free Paths

Our reduction makes extensive use of subdivisions, defined as follows.

**Definition 1.** The *complete subdivision* of a digraph $D$ is the digraph $S$ obtained from $D$ by replacing each arc $a \to b$ of $D$ by a new vertex $x_{ab}$ and the two arcs $a \to x_{ab} \to b$.

A digraph is a *subdivision digraph* if it is (isomorphic to) the complete subdivision of some digraph. A *core* of a subdivision digraph $S$ is a set $C$ of vertices of $S$ such that there exists a digraph $D$ with $S$ the complete subdivision of $D$ and $C = V(D)$.

Figure 3(a) shows a simple example of a complete subdivision. By definition, each subdivision digraph has a core. And it is easy to see that directed even cycles are the only weakly connected subdivision digraphs with more than one core. But we shall make no use of the latter observation. All we need are some simple facts about those vertices that are created through the subdivision process.

**Definition 2.** A vertex $x$ of a digraph is called *thin* if

$$\delta^+(x) = \delta^-(x) = 1.$$

**Fact 1.** *Let $C$ be a core of a subdivision digraph $S$ and let $\bar{C} := V(S) \setminus C$. Then all vertices in $\bar{C}$ are thin. And for each arc $a \to b$ of $S$, one of $a$ and $b$ lies in $C$ and the other in $\bar{C}$.*

Since our digraphs have no multiple edges, two different vertices of a subdivision digraph cannot share a common in-neighbour and a common out-neighbour.

**Fact 2.** *A subdivision digraph contains no two vertices $a \neq b$ with*

$$I(a) \cap I(b) \neq \emptyset \quad \text{and} \quad O(a) \cap O(b) \neq \emptyset.$$

Finally, let us have another look at Figure 2 from the proof of Proposition 1. The $k$-th root we have constructed there, consisted mainly of long isolated paths. Such paths will turn out very useful, so let us give them a name.

**Definition 3.** A path $a_1 \to a_2 \to \cdots \to a_k$ in a digraph is called *free* if

$$O(a_i) = \{a_{i+1}\}, \quad I(a_{i+1}) = \{a_i\} \quad \text{for } 1 \leq i < k.$$

So any non-path arc incident to the vertices of a free path either enters the starting vertex of that path or leaves its ending vertex.

### The Reduction in Detail

We are now able to state the final version of our reduction in detail. For a given pair $D_1, D_2$ of digraphs, we construct a digraph $\mathcal{R}(D_1, D_2)$ as follows.

1. Make $k-2$ isomorphic copies $D_3, \ldots, D_k$ of $D_2$.

2. Extend each $D_i$, $1 \leq i \leq k$, by two new vertices $r_i, s_i$ introducing the arcs
$$r_i \to a \quad \text{and} \quad s_i \to a \to s_i \quad \text{for each } a \in V(D_i).$$
This yields $k$ "suspensions" $\hat{D}_1, \ldots, \hat{D}_k$ (see Figure 3(b)).

3. Form the complete subdivision $\bar{D}_i$ of $\hat{D}_i$ for $1 \leq i \leq k$.

4. Let $\mathcal{R} := \bar{D}_1 \dot{\cup} \bar{D}_2 \dot{\cup} \cdots \dot{\cup} \bar{D}_k$ be the disjoint union of the $\bar{D}_i$.

This construction can clearly be done in polynomial time. Also observe that each $\hat{D}_i$ in the construction is weakly connected and hence so are the $\bar{D}_i$. Therefore, $\bar{D}_1, \ldots, \bar{D}_k$ are exactly the weakly connected components of $\mathcal{R}(D_1, D_2)$. Further note that because of the vertices $s_i$, each vertex of $\mathcal{R}(D_1, D_2)$ has a positive out-degree, and the only vertices with in-degree zero are the $r_i$, which all lie in different components.

## 4 Finding the Isomorphism

For the above reduction $\mathcal{R}$ we will now be able to show the desired implication

$$\mathcal{R}(D_1, D_2) \text{ has a } k\text{-th root} \Rightarrow D_1 \text{ and } D_2 \text{ are isomorphic.}$$

Taking into account the remarks about the construction $\mathcal{R}$ at the end of the preceding section and observing that the two digraphs $D_1$ and $D_2$ are isomorphic if and only if their subdivided suspensions $\bar{D}_1$ and $\bar{D}_2$ are isomorphic, we see that what we essentially need is the following proposition.

**Proposition 2.** *Let $D$ be a subdivision digraph with exactly $k$ weakly connected components $D_1, \ldots, D_k$. Every vertex of $D$ have positive out-degree and each $D_i$, $1 \leq i \leq k$, have exactly one vertex with in-degree zero. If $D$ has a $k$-th root then all the $D_i$ are isomorphic.*

It turns out that a root of a digraph $D$ as in Proposition 2 is just the one we constructed in our proof of Proposition 1. To show this, we have to provide several lemmas. The prevailing concept in all of these is that of a free path. These paths are used to establish local isomorphisms between the components of $D$.

**Lemma 1.** *Let $D$ be a subdivision digraph with a $k$-th root $R$. Each vertex have positive out-degree in $D$ and the set $Q := \{q \in V(D) \mid \delta_D^-(q) = 0\}$ have size $k$. Then $R[Q]$, the induced subgraph on $Q$, is a path (of length $k-1$). Moreover, this path is free in $R$.*

*Proof.* First observe that no $q \in Q$ lies on a cycle in $R$ since that would induce $\delta_D^-(q) > 0$. Second, we have

$$I_R(Q) \subseteq Q \tag{2}$$

because for any $x \in I_R(Q) \setminus Q$ with some $q \in Q$, we have $I_D(q) \supseteq I_R^{k-1}(x)$. The left-hand side is empty, the right-hand side not. Hence, such an $x$ cannot exist. We further fix a core $C$ of $D$ and note that $Q \subseteq C$.

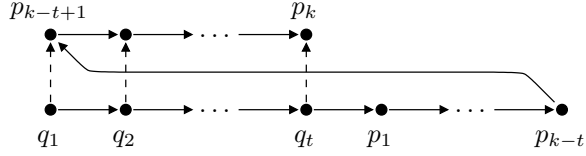

Figure 4: Finding the first free path. (The dashed arcs are from $D$.)

Now pick a longest path

$$q_1 \to q_2 \to \cdots \to q_t \tag{3}$$

in $R[Q]$. We shall see that $t = k$.

Note that

$$\delta_R^-(q_1) = 0 \tag{4}$$

because of (2), the maximality of the path (3), and the fact that there are no cycles in $R[Q]$.

Because each vertex has positive out-degree in $D$, there is a path

$$q_t = p_0 \to p_1 \to \cdots \to p_k \in O_D(q_t) \tag{5}$$

in $R$. We show

$$p_i \notin Q \quad \text{for } 1 \leq i \leq k \tag{6}$$

by induction. Clearly $p_1 \notin Q$ because the converse would contradict the maximality of the path (3) or yield a cycle in $Q$. And for $1 \leq i < k$, we have the chain of implications

$$p_i \notin Q \Rightarrow \delta_D^-(p_i) > 0 \Rightarrow \delta_D^-(p_{i+1}) > 0 \Rightarrow p_{i+1} \notin Q.$$

Combining the paths (3) and (5), we get $p_{k-t+1} \in O_D(q_1)$. Figure 4 shows a complete picture of the situation.

Because of $q_1 \in C$ and Fact 1, the vertex $p_{k-t+1}$ is thin and we conclude $I_D(p_{k-t+1}) = \{q_1\}$. This implies $I_R^{k-1}(p_{k-t}) = \{q_1\}$ and by (4) even $\delta_D^-(p_{k-t}) = 0$. So $p_{k-t} \in Q$, which, by (6), means $k = t$.

It remains to show that our path (3) is free. First, assume for contradiction that there is an arc $q_i \to v$ in $R$ with $v \neq q_{i+1}$, $1 \leq i < k$. We cannot have $v = q_0$ since then $q_0$ would lie on a cycle. So $v \neq q_0, q_{i+1}$. We claim $\delta_R^-(v) \geq 2$, i.e., $I_R(v) \supsetneq \{q_i\}$. If $v \in Q$, the path (3) yields a second in-neighbour of $v$. If $v \notin Q$, note that $\delta_D^-(q_{i+1}) = 0$ implies $I_R^{k-1}(q_i) = \emptyset$. So $I_D(v) = \{q_i\}$ would yield the contradiction $\delta_D^-(v) = 0$. Thus we see that actually $\delta_R^-(v) \geq 2$. Now consider any $z \in O_R^{k-1}(v)$. By Fact 1, $z$ is thin; a contradiction to $\delta_D^-(z) \geq \delta_R^-(v) \geq 2$.

Finally we check that there is no arc $u \to q_i$ in $R$ with $u \neq q_{i-1}$, $2 \leq i \leq k$. But similarly to the preceding argument, such an arc would give $q_{i-1}, u \in I_D(z)$ for any $z \in O_R^{k-1}(q_i)$; again a contradiction to $z$ thin. □

Lemma 1 will serve as a kind of induction basis: it gives us a first free path. The next four lemmas will allow us to propagate free paths through our subdivision digraph $\mathcal{R}(D_1, D_2)$.

**Lemma 2.** *Let $D$ be a subdivision digraph with a $k$-th root $R$ and let $C$ be a core of $D$. Let $R$ contain a free path $a_1 \to a_2 \to \cdots \to a_k$ with all $a_i \in C$, $1 \leq i \leq k$. Further assume that*

$$I_R^{k-1}(x) = \emptyset \Rightarrow O_R(x) \subseteq C \tag{7}$$

*for all $x \in V(D)$.*

*Then the sets $O_D(a_i)$ have the same cardinality $m$ for all $i \in \{1, \ldots, k\}$ and there exist $m$ disjoint free paths*

$$u_1^\mu \to u_2^\mu \to \cdots \to u_k^\mu, \quad 1 \leq \mu \leq m$$

*in $R$ with*

$$\bigcup_{\mu=1}^m u_i^\mu = O_D(a_i), \quad 1 \leq i \leq k.$$

Lemma 2 essentially states that locally, the root $R$ looks just as the one we constructed in Figure 2.

*Proof.* We start with the simple observation that all vertices in the set $\bigcup_{i=1}^k O_D(a_i)$ are thin because they neighbour vertices from $C$.

The proof consists of two steps. We first investigate $R$-*out*-neighbours of vertices in $N$, then $R$-*in*-neighbours.

*Claim* 1. Let $1 \leq i < k$. For each $x \in O_D(a_i)$, we have $O_R(x) = \{y\}$ with some $y \in O_D(a_{i+1})$.

Clearly $\delta_R^+(x) > 0$ because $x$ is thin. Pick some $y \in O_R(x)$. Since $O_R(a_i) = \{a_{i+1}\}$, we have $x \in O_R^{k-1}(a_{i+1})$ and thus $y \in O_D(a_{i+1})$ as claimed.

Now pick another vertex $z \in O_R(x)$; note that also $z \in O_D(a_{i+1})$. We have $O_R^{k-1}(y), O_R^{k-1}(z) \neq \emptyset$ because $y$ and $z$ are thin, and there exists a unique vertex $b$ with $O_D(x) = \{b\}$ because $x$ is also thin. So $y, z \in O_R(x)$ implies $O_R^{k-1}(y) = \{b\} = O_R^{k-1}(z)$ and thus $O_D(y) = O_D(z)$. Together with $y, z \in O_D(a_{i+1})$ and Fact 2, this yields $y = z$ and Claim 1 is proved.

A similar result can be obtained for the predecessors of the $a_i$.

*Claim* 2. Let $1 < i \leq k$. For each $x \in O_D(a_i)$ we have $I_R(x) = \{y\}$ with some $y \in O_D(a_{i-1})$.

Clearly $\delta_R^-(x) > 0$ because $x$ is thin. As above, we pick some $y \in I_R(x)$. We have $I_R^{k-1}(y) \subseteq I_R^k(x) = \{a_i\}$. Since $I_R^{k-1}(y) = \emptyset$ would, together with $x \notin C$, yield a contradiction to our precondition (7), we conclude $I_R^{k-1}(y) = \{a_i\}$ and thus $y \in O_D(a_{i-1})$.

We pick another vertex $z \in I_R(x)$. We have $O_R^{k-1}(x) \subseteq O_R^k(y) \cap O_R^k(z)$, and since the former set is nonempty, neither is the latter. Together with $y, z \in O_D(a_{i-1})$ and Fact 2, this yields $y = z$ and Claim 2 is proved.

Combining both claims, we get

$$O_R(O_D(a_i)) = O_D(a_{i+1}) \quad \text{for } 1 \leq i < k,$$
$$I_R(O_D(a_i)) = O_D(a_{i-1}) \quad \text{for } 1 < i \leq k.$$

Also each vertex in $\bigcup_1^{k-1} O_D(a_i)$ has a unique out-neighbour in $R$ and each vertex in $\bigcup_2^k O_D(a_i)$ has a unique in-neighbour in $R$. So all the sets $O_D(a_i)$, $1 \leq i \leq k$, have the same size $m$, and we thus get the $m$ disjoint free paths as stated in the lemma. □

**Lemma 3.** *Let $D$ be a subdivision digraph with a $k$-th root $R$ and let $C$ be a core of $D$. Let $R$ contain a free path $a_1 \to a_2 \to \cdots \to a_k$ with all $a_i \in C$, $1 \leq i \leq k$. Further assume $O_R^{k-1}(x) > 0$ for all $x \in V(D)$.*

*Then the sets $I_D(a_i)$ have the same cardinality $m$ for all $i \in \{1, \ldots, k\}$ and there exist $m$ disjoint free paths*

$$u_1^\mu \to u_2^\mu \to \cdots \to u_k^\mu, \quad 1 \leq \mu \leq m$$

*in $R$ with*

$$\bigcup_{\mu=1}^m u_i^\mu = I_D(a_i), \quad 1 \leq i \leq k.$$

*Proof.* This lemma is merely the transpose of Lemma 2.

If we reverse all arcs in $D$ and $R$ yielding digraphs $D'$ and $R'$, then in-neighbours become out-neighbours and vice versa. Clearly, $D'$ is still a subdivision digraph, $R'$ is a $k$-th root of $D'$, and $C$ is a core of $D'$. Free paths in $R$ have their orientation changed in $R'$ but they surely remain free. So this lemma follows from Lemma 2 if we also note that the precondition $O_R^{k-1}(x) > 0$ is even stronger than its counterpart (7). □

Lemmas 2 and 3 provided the "induction step" from a free path in a core $C$ to free paths through its $D$-neighbourhood, which of course lies in the complement of $C$. The next two lemmas take care of the other direction, from $C$'s complement to $C$.

**Lemma 4.** *Let $D$ be a subdivision digraph with a $k$-th root $R$. Let $C$ be a core of $D$ and let $R$ contain a free path $u_1 \to u_2 \to \cdots \to u_k$ with no $u_i$ lying in $C$. Let $a_i$ denote the unique $D$-out-neighbour of $u_i$ and assume that $\delta_D^+(a_i) > 0$ for each $i \in \{1, \ldots, k\}$. Then $R$ contains the free path $a_1 \to a_2 \to \cdots \to a_k$.*

*Proof.* Let $1 \leq i < k$. From $O_D(u_i) = \{a_i\}$ and $O_R(u_i) = \{u_{i+1}\}$ follows $O_R^{k-1}(u_{i+1}) = \{a_i\}$. This yields $O_R(a_i) = O_D(u_{i+1}) = \{a_{i+1}\}$.

It remains to show that $\delta_R^-(a_i) = 1$ for $1 < i \leq k$. Since $\delta_D^+(a_i) > 0$, there exists some $z \in O_R^{k-1}(a_i)$ and as a $D$-out-neighbour of $a_{i-1}$ such a $z$ is thin. We get $1 = \delta_D^-(z) \geq \delta_R^-(a_i) \geq 1$. Hence $I_R(a_i) = \{a_{i-1}\}$ as claimed. □

Reversing all arcs in Lemma 4, we immediately get the transpose statement with the roles of in- and out-neighbours interchanged.

**Lemma 5.** *Let $D$ be a subdivision digraph with a $k$-th root $R$. Let $C$ be a core of $D$ and let $R$ contain a free path $u_1 \to u_2 \to \cdots \to u_k$ with no $u_i$ lying in $C$. Let $a_i$ denote the unique $D$-in-neighbour of $u_i$ and assume that $\delta_D^-(a_i) > 0$ for each $i \in \{1, \ldots, k\}$. Then $R$ contains the free path $a_1 \to a_2 \to \cdots \to a_k$.* □

Lemmas 1 to 5 can now be used to show Proposition 2.

*Proof of Proposition 2.* First we pick some core $C$ of $D$ and observe that the $k$ vertices with zero in-degree must all lie in $C$. Lemma 1 yields a free path $p$ through these vertices.

Starting from this free path, we now repeatedly apply lemmas 2 to 5. Note that the precondition (7) of Lemma 2 is satisfied by $D$ since the only vertices $x$ with $I_R^{k-1}(x) = \emptyset$ can be the first $k-1$ vertices on the path $p$. And because $p$ is free and lies completely in $C$, all their $R$-out-neighbours lie also on $p$ and thus in $C$.

Each component $D_i$ is weakly connected and contains some vertex of $p$, so we will eventually cover the whole digraph $D$ with free paths. Again by lemmas 2 to 5, each of these paths touches each component $D_i$ exactly once and all do so in the same order. Therefore, all these free paths are disjoint.

Now recall that in all four lemmas 2 to 5, the free paths found respect the adjacencies of the $D_i$. So they actually establish isomorphisms between all $D_i$ as claimed. □

The proof of Theorem 1 is almost completely contained in our side remarks along the way to Proposition 2 and might already be clear by now. Nonetheless, we shall give it here for completeness.

*Proof of Theorem 1.* We have to show that the reduction $\mathcal{R}$ works as intended. That is, the digraph $\mathcal{R}(D_1, D_2)$ has a $k$-th root if and only if the digraphs $D_1$ and $D_2$ are isomorphic.

We have already discussed the easy direction in Section 3: Proposition 1 tells us that an isomorphism between $D_1$ and $D_2$ implies a $k$-th root of $\mathcal{R}(D_1, D_2)$.

For the other direction, compare our remarks about the construction of $\mathcal{R}(D_1, D_2)$ at the very end of Section 3 with the preconditions of Proposition 2. We see that the digraph $\mathcal{R}(D_1, D_2)$ meets all those requirements. Hence, a $k$-th root of $\mathcal{R}(D_1, D_2)$ implies an isomorphism between the subdivided suspensions $\bar{D}_1$ and $\bar{D}_2$ of $D_1$ respectively $D_2$. But as we have already observed before, $\bar{D}_1$ and $\bar{D}_2$ are isomorphic if and only if $D_1$ and $D_2$ are. $\square$

## 5 Towards Isomorphism Completeness

Reviewing Propositions 1 and 2, we immediately see that the $k$-th root problem is isomorphism complete for a certain subclass of directed graphs. Namely those digraphs that are created by our reduction.

**Theorem 2.** *For each $k \geq 2$, the problem of deciding whether a subdivision digraph with exactly $k$ weakly connected components whose vertices have all positive out-degree and such that each component contains exactly one vertex with in-degree zero, is isomorphism complete.*

*Proof.* Proposition 2 tells us that a subdivision digraph $D$ with exactly $k$ weakly connected components $D_1, \ldots, D_k$ and those additional conditions on the degrees has a $k$-th root if and only if all the $D_i$ are isomorphic. So we may reduce the $k$-th root problem for such digraphs to graph isomorphism as follows. Let $E_1$ be the disjoint union of $k - 1$ isomorphic copies of $D_1$ and let $E_2$ be the disjoint union of $D_2, \ldots, D_k$. Obviously, $E_1$ and $E_2$ are isomorphic if and only if all the $D_i$ are isomorphic. By Propositions 1 and 2, this is the case iff the subdivision digraph $D$ has a $k$-th root. $\square$

Admittedly, the preconditions of Theorem 2 are just too artificial to make it a remarkable fact of its own. We mention this byproduct rather as a possible starting point for further research. Maybe Theorem 2 can be generalized to more natural classes of directed graphs for which the $k$-th root problem remains isomorphism complete. What about, for example, complete subdivision graphs without any restrictions on connectedness and degrees?

## 6 Conclusions

We have shown that computing $k$-th roots of directed graphs (or equivalently, binary relations or Boolean matrices) is at least as hard as the graph isomorphism problem. The most important technique for this result are subdivisions. They essentially allowed us to derive uniqueness of a $k$-th root. The notion of free paths then led us from roots to local isomorphisms.

Since graph isomorphism is not known—nor very much expected—to be $\mathcal{NP}$-complete, it remains to find out the precise complexity of digraph roots in-between graph isomorphism and $\mathcal{NP}$-completeness.

We tried to mark a possible starting point for an attack on this question. Our hardness proof yielded, as a side effect, a class of digraphs for which root finding is even isomorphism complete. While being a rather technical statement, this byproduct might lead to larger, more natural classes of digraphs for which root finding could be isomorphism complete, even if the general problem should turn out $\mathcal{NP}$-complete.

## Acknowledgements

I want to thank Hein van der Holst and Mark de Longueville for helpful discussions.